\documentclass[a4paper,10pt]{article}
\usepackage{lineno,hyperref}
\usepackage{subfigure}
\usepackage{graphicx}
\usepackage{floatrow}
\usepackage{amsmath}
\usepackage{lipsum}
\usepackage{verbatim}
\usepackage{txfonts}
\usepackage{ulem}
\usepackage{tikz}
\usepackage{float}
\usepackage{xcolor}
\usepackage{multicol}
\usetikzlibrary{shapes,arrows}
\usepackage{multirow}
\usepackage{indentfirst}

\date{}

\title{\bf Chaotic dynamics of the sugarcane borer-two parasitoid agroecosystem with
seasonality}

\author{Marat Rafikov$^1$, Alexandre Molter$^2$, João Inácio Moreira Bezerra$^2$, Elvira Rafikova$^1$, Maria Cristina Varriale$^3$\\ \small{$^1$Center for Engineering, Modeling and Applied Social Sciences, Federal University of the ABC, Brazil}\\ \small{$^2$Post Graduate Program in Mathematical Modeling, Federal University of Pelotas, Brazil}\\  \small{$^3$Post Graduate Program in Applied Mathematics, Federal University of Rio Grande do Sul, Brazil}}

\begin{document}
\maketitle

\begin{abstract}
 \hspace{-0.5cm}Sugarcane production is a significant and profitable agribusiness sector in many countries. Nevertheless, this industry suffers significant losses from the sugarcane pests, among which the most important one is the sugarcane borer ({\it Diatraea saccharalis}).  This pest population is hard to be controlled due to its different life stages, thus biological control ( with more than one predator species)  can be applied. Therefore, in this work, we present and analyze a mathematical model that describes the dynamics of the sugarcane borer and its two different life stages parasitoids: eggs ({\it Trichogramma galloi}) and larval ({\it Cotesia flavipes}). First, a host-parasitoid model is used to obtain the population dynamics,  which also considers the influence of seasonal variations. Then, system simulations and bifurcation diagrams show that the introduction of seasonality perturbations causes complex dynamics and results in limit cycles and strange attractors. 
\end{abstract}

\hspace{0.5cm}{\bf keywords}:
Sugarcane borer; Parasitoids; Dynamics; Seasonality; Chaos


\section{Introduction}
Sugarcane is a global commodity produced throughout the tropics and used for sweeteners, biofuels, and a growing range of bioproducts (including bioplastics) \cite{WinNT}. Several pests affect the sugarcane's production, specially the sugarcane borer {\it Diatraea sacharalis} \cite{parraetal2002, postali2019applied}. Significant damages of the plant are caused by this peculiar insect, that lays eggs on the surface of the sugarcane leaves, while the larvae lives inside the sugarcane stalk, carving internal galleries. Because the larvae live hidden from the surface, the pesticide control becomes inefficient and biological control is the alternative for this case.

In biological control, the pest populations are reduced by the insertion of their natural enemies(predators, parasitoids and pathogens) in the environment. The sugarcane borer larvae population has been controlled since the 1970s by the parasitoid, \textit{Cotesia flavipes}. Recently, the \textit{Trichogramma galloi parasitoid} is a new alternative for biological control of its egg population \cite{parraetal2002, postali2019applied}. The works \cite{gamez2009observation, gamez2010open, venturino2007diseases, venturino2008biological} presented mathematical modeling in biological control applications of prey-predator and host-parasitoid.

In biological systems with more than two species, prey-predator interactions become complex and consequently harder to model prey invasions. Mathematical modeling can be an important feature providing information about the natural systems' stability \cite{goh2012management}, along with computational simulations, revealing the behavior of these complex systems and understanding how the prey interacts with other species in the environment. 

Mathematical models for biological control, with only one predator population, are considered in \cite{rafikov2012mathematical} and \cite{rafikov2014dynamical} for egg and larval parasitoid, respectively. The mathematical model of interactions between the sugarcane borer and its eggs and larval parasitoids was proposed by Rafikov and Silveira \cite{rafikov2013dynamics}. In Molnár et al. \cite{molnar2016two}, this model was used for the formulation of some scenarios of biological pest control. In both publications, the populations are described by four population compartments: sugarcane borer's eggs
population, sugarcane borer's larvae population, \textit{Trichogramma galloi}(eggs parasitoid) population, \textit{Cotesia Flavipes}(larval parasitoid) population.
 
As pointed out by \cite{rinaldi1993multiple}, periodic external forces are of great importance in ecological systems since environments of the population communities vary periodically. There are many systems that have a very simple dynamic behavior in the constant parameter case but become very complex (multiplicity of attractors, catastrophes, and chaos) when they are periodically perturbed.

Seasonality in biological systems has already been addressed in several works, such as \cite{rinaldi1993multiple, gakkhar2003chaos, gakkhar2003seasonally, altizer2006seasonality, zhang2012complexity, white2020seasonality, bezerra2021biological, stollenwerk2017hopf}. Meanwhile, the inclusion of two parasitoids (eggs and larval) populations, as well as seasonal variations in the sugarcane borer agroecosystem dynamics are novelties. Bezerra et al. in \cite{bezerra2021biological}, models sugarcane borer and its larval parasitoid \textit{(Cotesia flavipes)} interaction, considering the influence of the seasonal variations on the dynamics of the system. Their results show that this variation generates chaotic dynamics in the system.

The system model in the present paper is an extension of the mathematical model in \cite{rafikov2013dynamics}, with the addition of the parasitized egg and larvae population of the sugarcane borer. This addition improves the estimation and observation of the system parameters, as it is much easier to monitor sugarcane borer's parasitized eggs and parasitized larvae than adult parasitoid population in real conditions. Moreover, seasonality variations are introduced into the system dynamics resulting in chaotic behavior.

The paper is organized as follows. In Section 2, our six-dimensional, continuous-time, dynamical system is proposed. Section 3 is dedicated to the study of the system's local and global dynamics. Seasonal dynamics of the sugarcane borer parasitoid agroecosystem is considered in Section 4. Section 5 discusses the results of previous Sections and concludes this paper.


\section{Mathematical Model}

The proposed continuous-time mathematical model describes the interactions between the sugarcane borer, its egg and larval parasitoid, considering six population densities:
the un-parasitized egg population density of the sugarcane borer, $x_1$; the parasitized egg population density of the sugarcane borer,, $x_2$;  the density of the adult egg parasitoid \textit{Trichogramma galloi}, $x_3$; the unparasitized larvae density of the sugarcane borer, $x_4$; the parasitized larvae density of the sugarcane borer, $x_5$; the density of the adult larval parasitoid \textit{Cotesia flavipes}, $x_6$. So, the mathematical model has the following form:
\begin{align}
    \label{model}
    \frac{dx_1}{dt} &= rx_1\left(1 - \frac{x_1}{K}\right) - m_1x_1 - n_1x_1 - \alpha x_1x_3, \nonumber \\
    \frac{dx_2}{dt} &= \alpha x_1x_3 - m_2x_2 - n_2x_2, \nonumber \\
    \frac{dx_3}{dt} &= \gamma_1 n_2x_2 - m_3x_3, \nonumber \\
    \frac{dx_4}{dt} &= n_1x_1 - m_4x_4 - n_3x_4 - \beta x_4x_6, \\
    \frac{dx_5}{dt} &= \beta x_4x_6 - m_5x_5 - n_4x_5, \nonumber \\
    \frac{dx_6}{dt} &= \gamma_2n_4x_5 - m_6x_6. \nonumber
\end{align}

In the system of differential equations \eqref{model}, the 16 parameters are defined as follows. $r$ is the intrinsic oviposition rate of female sugarcane borer; $K$ is the potential
maximum of oviposition rate of female sugarcane borer; $m_1, m_2, m_3, m_4, m_5$ and $m_6$ are the mortality rates of the un-parasitized egg, parasitized egg, egg parasitoid, un-
parasitized larvae, parasitized larvae and larvae parasitoid populations, respectively; $n_1$ is the fraction of the sugarcane borer larvae population which emerges from the eggs per unit of time; $n_2$ is the fraction of the parasitized egg population from which larval
parasitoids emerge in a time unit; $n_3$ is the fraction of the un-parasitized sugarcane borer
larvae from which pupae emerge in a time unit; $n_4$ is the fraction of the parasitized
sugarcane borer larvae from which larvae parasitoids emerge in a time unit; $\alpha$ and $\beta$ are are
the intrinsic parasitism rate of the egg and larvae parasitoids, respectively; $\gamma_1$ and $\gamma_2$ are
numbers of adult parasitoids which emerge from a unit of parasitized eggs and larvae,
respectively.

\section{System Equilibrium states and their stability}

\subsection{Equilibrium states}

To obtain the equilibrium points, we equal the right-hand sides of the system \eqref{model} to zero. We obtain five following equilibrium by sequences of ${x_i, i = 1, 2, \cdots, 6}$.
\begin{itemize}
    \item Extinction of all populations: $E_1 = (0, 0, 0, 0, 0, 0)$.
    \item Extinction of all parasitoid and parasitized populations: $E_2 = \left(\frac{K}{r}(r - m_1 - n_1), 0, 0, \frac{Kn_1(r - m_1 - n_1)}{r(m_4+n_3)}, 0, 0\right)$.
    \item Extinction of the egg parasitoid and parasitized egg populations: $E_3 = \left(\frac{K}{r}(r - m_1 - n_1), 0, 0, p_4^{\ast}, p_5^{\ast}, p_6^{\ast}\right)$.
    \item Extinction of the larvae parasitoid and parasitized larvae populations: $E_4 = (x_1^{\ast}, x_2^{\ast}, x_3^{\ast}, q^{\ast}, 0, 0)$.
    \item Coexistence of all populations: $E_5 = (x_1^{\ast}, x_2^{\ast}, x_3^{\ast}, x_4^{\ast}, x_5^{\ast}, x_6^{\ast}).$
\end{itemize}
The values cited above are given as follows.
\begin{align*}
    x_1^{\ast} &= \frac{m_3(m_2 + n_2)}{\alpha \gamma_1 n_2}, & x_{3}^{\ast} &= \frac{1}{\alpha} \left[r\left(1 - \frac{m_3(m_2 + n_2)}{\alpha \gamma_1 n_2 K}\right) - m_1 - n_1\right], & x_2^{\ast} &= \frac{m_3}{\gamma_1 n_2}x_3^{\ast}, \\
    x_4^{\ast} &= \frac{m_6(m_5 + n_4)}{\beta \gamma_2 n_4}, & x_5^{\ast} &= \frac{m_3n_1(m_2 + n_2)}{\alpha \gamma_1 n_2(m_5 + n_4)} - \frac{m_6(m_4 + n_3)}{\beta \gamma_2 n_4}, & x_6^{\ast} &= \frac{\gamma_2n_4}{m_6}x_5^{\ast}, &  \\
    p_4^{\ast} &= \frac{m_6(m_5 + n_4)}{\beta \gamma_2 n_4}, & p_5^{\ast} &= \frac{n_1}{m_5 + n_4}\left[\frac{K}{r}(r - m_1 - n_1)\right] - \frac{m_6(m_4 + n_3)}{\beta \gamma_2 n_4}, & p_6^{\ast} &= \frac{\gamma_2n_4}{m_6}p_5^{\ast},\\
    q^{\ast} &= \frac{n_1}{m_4 + n_3}x_1^{\ast}.
\end{align*}

Since we are modelling a biological system, where the dependent variables are populations, the space phase includes only positive or zero values for all coordinates, thus defining their biological
viability. When analyzing, in next section, each one of these equilibrium points, conditions will be set for their biological viability.

\subsection{Local stability analysis of equilibrium points}

In order to study the local the stability of these equilibrium states, system \eqref{model} is linearized in a small neighborhood of each equilibrium state, and the Jacobian matrix being computed as:
\begin{equation}
    \label{jacobian}
    J=
    \begin{bmatrix} 
	a_{11} & 0 & a_{13} & 0 & 0 & 0\\
	a_{21} & a_{22} & a_{23} & 0 & 0 & 0\\
	0 & a_{32} & a_{33} & 0 & 0 & 0\\
	a_{41} & 0 & 0 & a_{44} & 0 & a_{46}\\
	0 & 0 & 0 & a_{54} & a_{55} & a_{56}\\
	0 & 0 & 0 & 0 & a_{65} & a_{66}\\
	\end{bmatrix},
\end{equation}
in which:
\begin{align*}
    a_{11} &= r - \frac{2rx_1}{K} - m_1 - n_1 - \alpha x_3, a_{13} = -\alpha x_1, a_{21} = \alpha x_3, a_{22} = -m_2 - n_2, \\
    a_{23} &= \alpha x_1, a_{32} = \gamma_1n_2, a_{33} = -m_3, a_{41} = n_1, a_{44} = -m_4 - n_3 - \beta x_6, \\
    a_{46} &= -\beta x_4, a_{54} = \beta x_6, a_{55} = -m_5 - n_4, a_{56} = \beta x_4, a_{65} = \gamma_2n_4, a_{66} = -m_6.
\end{align*}

The matrix \eqref{jacobian} can be written as a block matrix:
\begin{equation}
    \label{jacobian_block}
    \begin{bmatrix} 
	A & 0\\
	C & B \\
	\end{bmatrix},
\end{equation}
in which:
\begin{align*}
    A &= \begin{bmatrix}
    a_{11} & 0 & a_{13}\\
	a_{21} & a_{22} & a_{23}\\
	0 & a_{32} & a_{33}\\
    \end{bmatrix},
    C = 
    \begin{bmatrix}
     a_{41} & 0 & 0 \\
     0 & 0 & 0 \\
     0 & 0 & 0 \\
    \end{bmatrix}, 
    B = 
    \begin{bmatrix}
     a_{44} & 0 & a_{46} \\
     a_{54} & a_{55} & a_{56} \\
     0 & a_{65} & a_{66} \\
    \end{bmatrix},
\end{align*}
and 0 is a matrix with all elements equal to zero.

Matrices of the form \eqref{jacobian_block} are called block-lower-triangular. The full explanation on the block-triangular matrix determinant computation method can be found in \cite{gantmachertheory}, and  briefly described below: 

If D is a block-triangular matrix, then the determinant of the matrix is equal to the product of determinant of diagonal cells:
\begin{equation}
\label{detD}
    \det{(D)} = \det{(D_{11})}\det{(D_{22})}\cdots\det{(D_{nn})}. 
\end{equation}
The matrix $D = J - \lambda I$, where $I$ is the identity matrix of dimensions $n \times n$, is block-lower-triangular too. Using the above mentioned rule, the characteristic equation is given as follows:
\begin{equation}
    \label{char_eq}
    \det{(D)} = \det{(A - \lambda I)}\det{(B - \lambda I)} = 0.
\end{equation}
Hence, the characteristic equation \eqref{char_eq} is given by:
\begin{equation}
    \label{char_eq_expanded}
    \begin{vmatrix}
    a_{11} - \lambda & 0 & a_{13}\\
	a_{21} & a_{22} - \lambda & a_{23}\\
	0 & a_{32} & a_{33} -\lambda \\
    \end{vmatrix}
    \begin{vmatrix}
    a_{44}-\lambda & 0 & a_{46} \\
     a_{54} & a_{55}-\lambda & a_{56} \\
     0 & a_{65} & a_{66}-\lambda \\
    \end{vmatrix}
    = 0.
\end{equation}
From \eqref{char_eq_expanded}, we get that:
\begin{equation}
    \label{char_1}
    \begin{vmatrix}
    a_{11} - \lambda & 0 & a_{13}\\
	a_{21} & a_{22} - \lambda & a_{23}\\
	0 & a_{32} & a_{33} -\lambda \\
    \end{vmatrix} 
    = 0,
\end{equation}

\begin{equation}
    \label{char_2}
    \begin{vmatrix}
    a_{44}-\lambda & 0 & a_{46} \\
     a_{54} & a_{55}-\lambda & a_{56} \\
     0 & a_{65} & a_{66}-\lambda \\
    \end{vmatrix}
    = 0.
\end{equation}

Now, applying this rule to the equilibrium points found above, we obtain for the equilibrium point $E_1 = (0, 0, 0, 0, 0, 0)$, the matrix \eqref{jacobian} has a triangular form, and the eigenvalues are given by:
$\lambda_1= r - m_1 - n_1,\,\lambda_2 = -m_2 - n_2,\, \lambda_3 = -m_3,\, \lambda_4 = -m_4 - n_3,\, \lambda_5 = -m_5 - n_4,\,\text{and}\,\lambda_6= -m_6.$
Therefore, it follows that equilibrium $E_1$ is asymptotically stable if $r < m_1 + n_1$.

For the equilibrium $E_2 = \left(\frac{K}{r}(r - m_1 - n_1), 0, 0, \frac{Kn_1(r - m_1 - n_1)}{r(m_4+n_3)}, 0, 0\right)$, which is biologically viable if $r > m_1 + n_1$, the parameters of the characteristic equation \eqref{char_eq} are given as:
\begin{align*}
    a_{11} &= -(r - m_1 - n_1), a_{13} = -\frac{aK}{r}(r - m_1 - n_1), a_{21} = 0, a_{22} = -m_2 - n_2, \\
    a_{23} &= \frac{aK}{r}(r - m_1 - n_1), a_{32} = \gamma_1n_2, a_{33} = -m_3, a_{44} = -m_4 - n_3, \\
    a_{46} &= -\frac{\beta Kn_1(r - m_1 - n_1)}{r(m_4 + n_3)}, a_{54} = 0, a_{55} = -m_5 - n_4, a_{56} = \frac{\beta Kn_1(r - m_1 - n_1)}{r(m_4 + n_3)}, \\
    a_{65} &= \gamma_2n_4, a_{66} = -m_6.
\end{align*}

From \eqref{char_1} and \eqref{char_2}, we obtain:
\begin{align}
    \label{cond_1}
    (a_{11} - \lambda)[\lambda^2 - (a_{22} + a_{33})\lambda + a_{22}a_{33} - a_{23}a_{32}] &= 0 \\
    \label{cond_2}
    (a_{44} - \lambda)[\lambda^2 - (a_{55} + a_{66})\lambda + a_{55}a_{66}  - a_{56}a_{65}] &= 0
\end{align}

The Routh-Hurwitz criterion states that the eigenvalues of the second-degree
polynomial have negative real parts if, and only if, both coefficients are positive.
Analyzing \eqref{char_1} and \eqref{char_2}, we can conclude that when 

(a) $a_{11} < 0$,\;\;\; (b) $a_{22} + a_{33} - a_{23}a_{32} > 0$,\;\;\; (c) $a_{55} + a_{66} - a_{56}a_{65} > 0$,

\noindent then all the eigenvalues of the equation \eqref{char_eq_expanded} have negative real parts.

From conditions (a), (b) and (c), we obtain:
\begin{align}
    \label{eig_2}
    \alpha &< \frac{rm_3(m_2+n_2)}{\gamma_1n_2K(r - m_1 - n_1)} \\
    \label{eig_3}
    \beta &< \frac{rm_6(m_4 + n_3)(m_5 + n_4)}{\gamma_2n_1n_4K(r - m_1 - n_1)} \\
    \label{eigr}
    r &> m_1 + n_1
\end{align}

Similarly, the equilibrium point $E_2$ is asymptotically stable if, and only if, the inequalities \eqref{eig_2}, \eqref{eig_3} and \eqref{eigr} are satisfied.

The equilibrium point $ E_3 = \left(\frac{K}{r}(r - m_1 - n_1), 0, 0, p_4^{\ast}, p_5^{\ast}, p_6^{\ast}\right)$ is biologically viable if:
\begin{align}
    \label{cond_e3}
    \beta &> \frac{rm_6(m_4 + n_3)(m_5 + n_4)}{\gamma_2n_1n_4K(r - m_1 - n_1)}, \\
    r &> m_1 + n_1. \nonumber 
\end{align}

The characteristic equation at $E_3$ can be written in the form \eqref{char_eq_expanded}, where:
\begin{align*}
    a_{11} &= -(r - m_1 - n_1), a_{13} = -\frac{aK}{r}(r - m_1 - n_1), a_{21} = 0, a_{22} = -m_2 - n_2, \\
    a_{23} &= \frac{aK}{r}(r - m_1 - n_1), a_{32} = \gamma_1n_2, a_{33} = -m_3, a_{44} = -m_4 - n_3 - \beta p_6^{\ast}, \\
    a_{46} &= -\beta p_4^{\ast}, a_{54} = \beta p_6^{\ast}, a_{55} = -m_5 - n_4, a_{56} = \beta p_4^{\ast}, a_{65} = \gamma_2n_4, a_{66} = -m_6.
\end{align*}

The parameter values of determinant \eqref{char_1} are the same as $E_2$, and the conditions \eqref{eig_2} and \eqref{eigr} are satisfied.

Considering the determinant \eqref{char_2}
\begin{equation}
    \label{char_2_eq3}
    \begin{vmatrix}
    a_{44}-\lambda & 0 & a_{46} \\
     a_{54} & a_{55}-\lambda & a_{56} \\
     0 & a_{65} & a_{66}-\lambda \\
    \end{vmatrix}
    = 0,
\end{equation}
we have
\begin{equation}
\label{char_poly_e3}
    \lambda^3 + b_1\lambda^2 + b_2\lambda + b_3 = 0,
\end{equation}
where
\begin{align}
\label{stab_e3}
    b_1 &= m_4 + n_3 + \beta p_6^{\ast} > 0, \;\;\;\;\;
    b_2= (m_4 + n_3 + \beta p_6^{\ast})(m_5 + n_4) > 0, \nonumber \\
    b_3 &= \beta m_6(m_5 + n_4)p_6^{\ast} > 0, \;\;\;\;\;
    b_1b_2 - b_3 > 0. 
\end{align}

Therefore, we obtain that the equilibrium point $E_3$ is asymptotically
stable if, and only if, the inequalities \eqref{eig_2}, \eqref{cond_e3} and \eqref{eigr} are satisfied.

The equilibrium point $E_4 = (x_1^{\ast}, x_2^{\ast}, x_3^{\ast}, q^{\ast}, 0, 0)$ is biologically viable if
\begin{align}
    \label{cond_e4}
    \alpha &> \frac{rm_3(m_2+n_2)}{\gamma_1n_2K(r - m_1 - n_1)}, \\
    r &> m_1+n_1. \nonumber 
\end{align}

The characteristic equation at $E_4$ can be written in the form \eqref{char_eq_expanded} where:
\begin{align}
\label{param_values_e4}
    a_{11} &= -\frac{rx_1^{\ast}}{K}, a_{13} = -ax_1^{\ast}, a_{21} = ax_3^{\ast}, a_{22} = -m_2 - n_2, \nonumber \\
    a_{23} &= ax_1^{\ast}, a_{32} = \gamma_1n_2, a_{33} = -m_3, a_{44} = -m_4 - n_3, \\
    a_{46} &= -\beta q^{\ast}, a_{54} = 0, a_{55} = -m_5 - n_4, a_{56} = \beta q^{\ast}, a_{65} = \gamma_2n_4, a_{66} = -m_6. \nonumber
\end{align}

Considering the determinant \eqref{char_1} with parameter values \eqref{param_values_e4} we have:
\begin{equation}
    \label{char_poly_e4}
    \lambda^3 + c_1\lambda^2 + c_2\lambda + c_3 = 0.
\end{equation}
where
\vspace{-0.2cm}
$$c_1= m_2 + m_3 + n_2 + \frac{rx_1^{\ast}}{K} > 0,\;
c_2=(m_2 + m_3 + n_2)\frac{rx_1^{\ast}}{K} > 0, \;
c_3=\alpha m_3(m_2 + n_2)x_3^{\ast} > 0.$$

From $c_1c_2 - c_3 > 0$, we obtain:
\begin{equation}
    \label{cond_c4}
    \alpha < \frac{rm_3(m_2 + n_2)}{\gamma_1n_2zK},
\end{equation}
where
\begin{align*}
    z &= -\frac{h_1}{2} + \sqrt{\frac{h_1^2}{4} + h_2}\,, \\
    h_1 &= m_2 + m_3 + n_2 + \frac{m_3(m_2 + n_2)}{m_2 + m_3 + n_2}, \\
    h_2 &= \frac{m_3(m_2 + n_2)(r - m_1 - n_1)}{m_2 + m_3 + n_2}.
\end{align*}

Considering the determinant \eqref{char_2} with parameter values \eqref{param_values_e4} we have:
\begin{equation}
    \label{char_2_e4}
    (a_{44} - \lambda)[\lambda^2 - (a_{55} + a_{66})\lambda  + a_{55}a_{66} - a_{56}a_{65}] = 0,
\end{equation}
where \;
   $ a_{44} = -m_4 - n_3 < 0, -(a_{55} + a_{66}) > 0.$

From $a_{55}a_{66} - a_{56}a_{65} > 0$ we obtain:
\begin{equation}
    \label{beta_eq4_cond}
    \beta < \frac{m_6(m_5+n_4)}{\gamma_2n_4q^{\ast}} = \frac{\alpha\gamma_1n_2m_6(m_4+n_3)(m_5+n_4)}{\gamma_2n_1n_4m_3(m_2+n_2)}.
\end{equation}

The equilibrium point $E_4$ is asymptotically stable if, and only if, the following
inequalities are satisfied:
\begin{equation}
    \label{cond1_e4}
    \frac{rm_3(m_2 + n_2)}{\gamma_1n_2K(r - n_1 - m_1)} < \alpha < \frac{rm_3(m_2 + n_2)}{\gamma_1n_2zK},
\end{equation}

\begin{equation}
    \label{cond2_e4}
    \beta < \frac{\alpha\gamma_1n_2m_6(m_4+n_3)(m_5+n_4)}{\gamma_2n_1n_4m_3(m_2+n_2)}.
\end{equation}

Consider the equilibrium point $E_5 = (x_1^{\ast}, x_2^{\ast}, x_3^{\ast}, x_4^{\ast}, x_5^{\ast}, x_6^{\ast})$ where
\begin{align}
\label{eq5coor}
     x_1^{\ast} &= \frac{m_3(m_2 + n_2)}{\alpha \gamma_1 n_2}, x_{3}^{\ast} = \frac{1}{\alpha} \left[r\left(1 - \frac{m_3(m_2 + n_2)}{\alpha \gamma_1 n_2 K} - m_1 - n_1\right)\right], x_2^{\ast} = \frac{m_3}{\gamma_1 x_2}x_3^{\ast}, \\
    x_4^{\ast} &= \frac{m_6(m_5 + n_4)}{\beta \gamma_2 n_4}, x_5^{\ast} = \frac{m_3n_1(m_2 + n_2)}{\alpha \gamma_1 n_2(m_5 + n_4)} - \frac{m_6(m_4 + n_3)}{\beta \gamma_2 n_4}, x_6^{\ast} = \frac{\gamma_2n_4}{m_6}x_5^{\ast} \nonumber. 
\end{align}

The equilibrium point $E_5$ is biologically viable if $x_3^{\ast} > 0$ and $x_5^{\ast} > 0$. So, we get:
\begin{align}
\label{equilibrio5}
    \alpha &> \frac{rm_3(m_2 + n_2)}{\gamma_1n_2K(r - n_1 - m_1)}, \nonumber \\
    \beta &> \frac{\alpha\gamma_1n_2m_6(m_4+n_3)(m_5+n_4)}{\gamma_2n_1n_4m_3(m_2+n_2)}, \\
    r &> m_1 + n_1. \nonumber
\end{align}

The characteristic equation at $E_5$ can be written in the form \eqref{char_eq_expanded} where:
\begin{align}
\label{e5eq}
    a_{11} &= -\frac{rx_1^{\ast}}{K}, a_{13} = -ax_1^{\ast}, a_{21} = ax_3^{\ast}, a_{22} = -m_2 - n_2, \nonumber \\
    a_{23} &= ax_1^{\ast}, a_{32} = \gamma_1n_2, a_{33} = -m_3, a_{44} = -m_4 - n_3 - \beta x_6^{\ast}, \\
    a_{46} &= -\beta x_4^{\ast}, a_{54} = \beta x_6^{\ast}, a_{55} = -m_5 - n_4, a_{56} = \beta x_4^{\ast}, a_{65} = \gamma_2n_4, a_{66} = -m_6 \nonumber.
\end{align}

The parameter values of determinant \eqref{char_1} are the same of $E_4$, and the inequality \eqref{cond_c4} is satisfied.

Considering the determinant \eqref{char_2}, we obtain:
\begin{equation}
\label{char_poly_e5}
    \lambda^3 + g_1\lambda^2 + g_2\lambda + g_3 = 0,
\end{equation}
where
\begin{align}
    \label{poly_e5_coefficients}
    g_1 &= m_4 + m_5 + m_6 + n_3 + n_4 + \beta x_6^{\ast} > 0,\nonumber \\
    g_2 &=(m_4 + n_3 + \beta x_6^{\ast})(m_5 + n_4 + m_6) > 0, \\
    g_3 &= \beta m_6(m_5 + n_4)x_6^{\ast} > 0,\;\;\;
    g_1g_2 - g_3 > 0. \nonumber
\end{align}

Therefore, we obtain that equilibrium point $E_5$ is asymptotically
stable if, and only if, the following inequalities are satisfied:
\begin{equation}
    \label{cond1_e5}
    \frac{rm_3(m_2 + n_2)}{\gamma_1n_2K(r - m_1 - n_1)} < \alpha < \frac{rm_3(m_2 + n_2)}{\gamma_1n_2zK},
\end{equation}

\begin{equation}
    \label{cond2_e5}
    \beta > \frac{\alpha\gamma_1n_2m_6(m_4+n_3)(m_5+n_4)}{\gamma_2n_1n_4m_3(m_2+n_2)}.
\end{equation}

We can now summarize the results of this local stability analysis, after defining the following dimensionless parameters related to our resulting conditions (regions in the parameter space) for biological viability (b. v.) and for local stability (l.s.) of each equilibrium point:

\begin{align}
\label{eqmc}
    A_1 &\equiv \frac{r}{m_1 + n_1};\ A_2 \equiv \frac{\alpha \gamma_1n_2K(r - m_1 - n_1)}{rm_3(m_2+n_2)};\ A_3 \equiv \frac{\beta \gamma_2n_1n_4K(r-m_1-n_1)}{rm_6(m_4+n_3)(m_5+n_4) }; \\
    A_4 &\equiv \frac{\beta\gamma_2n_1n_4m_3(m_2+n_2)}{\alpha\gamma_1n_2m_6(m_4+n_3)(m_5+n_4)} = \frac{A_3}{A_2};\ A_5 \equiv \frac{\alpha\gamma_1n_2zK}{rm_3(m_2+n_2)} = \frac{A_2z}{r-m_1-n_1} \nonumber.
\end{align}

With these dimensionless parameters, whose critical value 1 is associated with a bifurcation in the behavior of the system, all conditions
previously deduced can be presented in a summary and complete form as shown in Table \ref{tab:my-table}.

\begin{table}[H]
\scriptsize
\begin{tabular}{llll|ll|ll|ll|ll|ll|}
\cline{5-14}
 &
   &
   &
   &
  \multicolumn{2}{c|}{$E_1$} &
  \multicolumn{2}{c|}{$E_2$} &
  \multicolumn{2}{c|}{$E_3$} &
  \multicolumn{2}{c|}{$E_4$} &
  \multicolumn{2}{c|}{$E_5$} \\ \hline
\multicolumn{4}{|c|}{$A_1 < 1$} &
  \multicolumn{1}{c|}{{b.v.}} &
  l.s. &
  \multicolumn{2}{c|}{{Not b.v.}} &
  \multicolumn{2}{c|}{{Not b.v.}} &
  \multicolumn{2}{c|}{{Not b.v.}} &
  \multicolumn{2}{c|}{{Not b.v.}} \\ \hline
\multicolumn{1}{|c|}{\multirow{5}{*}{\small $A_1 > 1$}} &
  \multicolumn{1}{c|}{\multirow{2}{*}{\small $A_2 < 1$}} &
  \multicolumn{2}{c|}{$A_3 < 1$} &
  \multicolumn{1}{c|}{b.v.} &
  un. &
  \multicolumn{1}{c|}{{b.v.}} &
  l.s. &
  \multicolumn{2}{c|}{{Not b.v.}} &
  \multicolumn{2}{c|}{{Not b.v.}} &
  \multicolumn{2}{c|}{{Not b.v.}} \\ \cline{3-14} 
\multicolumn{1}{|c|}{} &
  \multicolumn{1}{c|}{} &
  \multicolumn{2}{c|}{$A_3 > 1$} &
  \multicolumn{1}{c|}{b.v.} &
  un. &
  \multicolumn{1}{c|}{b.v.} &
  un. &
  \multicolumn{1}{c|}{{b.v.}} &
  l.s. &
  \multicolumn{2}{c|}{{Not b.v.}} &
  \multicolumn{2}{c|}{{Not b.v.}} \\ \cline{2-14} 
\multicolumn{1}{|c|}{} &
  \multicolumn{1}{c|}{\multirow{3}{*}{$A_2 > 1$}} &
  \multicolumn{1}{c|}{\multirow{2}{*}{$A_5 < 1$}} &
  $A_4<1$ &
  \multicolumn{1}{c|}{b.v.} &
  un. &
  \multicolumn{1}{c|}{b.v.} &
  un. &
  \multicolumn{1}{c|}{b.v.} &
  un. &
  \multicolumn{1}{c|}{{b.v.}} &
  l.s. &
  \multicolumn{2}{c|}{Not b.v.} \\ \cline{4-14} 
\multicolumn{1}{|c|}{} &
  \multicolumn{1}{c|}{} &
  \multicolumn{1}{c|}{} &
  $A_4>1$ &
  \multicolumn{1}{c|}{b.v.} &
  un. &
  \multicolumn{1}{c|}{b.v.} &
  un. &
  \multicolumn{1}{c|}{b.v.} &
  un. &
  \multicolumn{1}{c|}{b.v.} &
  un. &
  \multicolumn{1}{c|}{{b.v.}} &
  l.s. \\ \cline{3-14} 
\multicolumn{1}{|c|}{} &
  \multicolumn{1}{c|}{} &
  \multicolumn{1}{c|}{$A_5 > 1$} &
  $A_4>1$ &
  \multicolumn{1}{c|}{b.v.} &
  un. &
  \multicolumn{1}{c|}{b.v.} &
  un. &
  \multicolumn{1}{c|}{b.v.} &
  un. &
  \multicolumn{1}{c|}{b.v.} &
  un. &
  \multicolumn{1}{c|}{b.v.} &
  un. \\ \hline
\end{tabular}
\caption{Local stability and biological viability of the system equilibria, in which b.v. means biologically viable, l.s. means locally stable and un. means unstable.}
\label{tab:my-table}
\end{table}

\subsection{Global stability analysis of the coexistence equilibrium}

From local stability analysis, we get the equilibrium point $E_5$ is the only one with no null population density. Therefore this is the point of interest for the forthcoming global stability analysis.

We define a Lyapunov function as follows:
\begin{equation}
    \label{lyap_function}
    V(x_1,x_2,x_3,x_4,x_5,x_6) = \int_{x_1^{\ast}}^{x_1} \frac{y - x_1^{\ast}}{y} dy + \sum_{i = 2}^{6} \frac{(x_i - x_1^{\ast})^2}{2}.
\end{equation}


At the equilibrium point $E_5$, $V$ is zero, and it is positive for all other biologically viable equilibria.  The function is also radially unbounded, i.e., $V \rightarrow \infty$ when $x \rightarrow \infty$.

We can write the time derivative of $V$ along \eqref{model} as 
\begin{equation}
\label{derivative}
    \dot{V} = e^TPe, 
\end{equation}
where the matrix $P$ is\\
\begin{equation*}
    \begin{bmatrix} 
	-\frac{r}{K} & 0 & -\alpha & 0 & 0 & 0\\
	\alpha x_3 & -m_2-n_2 & \alpha x_1^\ast & 0 & 0 & 0\\
	0 & \gamma_1n_2 & -m_3 & 0 & 0 & 0\\
	n_1 & 0 & 0 & -m_4-n_3-\beta x_6 & 0 & -\beta x_4^\ast\\
	0 & 0 & 0 & \beta x_6 & -m_5 - n_4 & \beta x_4^\ast\\
	0 & 0 & 0 & 0 & \gamma_2n_4 & -m_6\\
	\end{bmatrix},
\end{equation*}
and the elements of vector $e$ are:
   $ e_i = x_i - x_i^{\ast}, i = 1, 2, ..., 6.$

If conditions \eqref{cond1_e5} and \eqref{cond2_e5} are satisfied then the matrix $P$ in \eqref{derivative} is negative definite, as is the time derivative of $V$ along the trajectories of \eqref{model}, and consequently the equilibrium point $E_5$ is globally asymptotically stable.

\begin{equation}
\alpha_c \equiv \frac{rm_3(m_2+n_2)}{\gamma_1n_2zK}.
    \label{alpha_c}
\end{equation}

\subsection{Hopf bifurcation analysis}

From Table \ref{tab:my-table}, it is evident the role of the dimensionless parameter $A_5$ in determining the region in the parameter space in which each one of the equilibrium points $E_4$ and $E_5$ is asymptotically stable. Furthermore, from the definition of $A_5$ in \eqref{eqmc}, it can be seen that the condition $A_5 < 1$ can be written equivalently in the form of $\alpha < \alpha_c$, where the critical value $\alpha_c$ is defined by
\begin{equation}
    \label{hopf_cond_1}
    \alpha_{c}  \equiv \frac{rm_3(m_2 + n_2)}{\gamma_1n_2zK}
\end{equation}

In \eqref{hopf_cond_1} we observe that the condition is exactly the same one as we got in \eqref{cond_c4}, from the characteristic \eqref{char_poly_e4}, for the local stability of the equilibria $E_4$ and $E_5$.

When $A_5 > 1$, that is, $\alpha > \alpha_c$, the positive coexistence equilibrium $E_5$ becomes unstable and a Hopf bifurcation occurs.

Now we can analyze the bifurcation of the model \eqref{model} assuming $\alpha$ as the bifurcation parameter and considering only three 
first equations of the system \eqref{model} which are not dependent on the variables $x_4, x_5$ and $x_6$. The traditional Hopf bifurcation criterion is stated in terms of the properties of the 
eigenvalues. Alternatively, Liu (1994) \cite{liu1994criterion} presented a criterion of Hopf bifurcation without
using the eigenvalues of the characteristic equation.  Liu’s approach is the one that is applied in the present Hopf bifurcation analysis, as follows:




\textit{Liu's criterion.} If the characteristic equation of the positive equilibrium point is given by:
   $ \lambda^3 + c_1(\alpha)\lambda^2 + c_2(\alpha)\lambda + c_3(\alpha) = 0,$
where $c_1(\alpha), c_2(\alpha)$ and $c_3(\alpha)$ are smooth functions of $\alpha$ in an open interval about $\alpha_c \in \mathcal{R}$ such that:
\begin{enumerate}
    \item [(a)] $c_1(\alpha_c) > 0, \Delta(\alpha_c) = c_1(\alpha_c)c_2(\alpha_c) - c_3(\alpha_c) = 0, c_3(\alpha_c) > 0$.
    \item [(b)] $\left(\frac{d\Delta}{d\alpha}\right)_{\alpha = \alpha_c} \neq 0$,
\end{enumerate}
then a simple Hopf bifurcation occurs at $\alpha = \alpha_c$.

Applying the Liu’s criterion to the characteristic equation \eqref{char_poly_e4}, we observe that
    $$c_1 = m_2 + m_3 + n_2 + \frac{rx_1^{\ast}}{K} > 0,\;\;
    c_2 = (m_2 + m_3 + n_2)\frac{rx_1^{\ast}}{K} > 0, \;\;
    c_3 = \alpha(m_2 + n_2)x_3^{\ast} > 0,$$
for all positive values of $\alpha$.

Solving the equation
 $c_1(\alpha_c)c_2(\alpha_c) - c_3(\alpha_c) = 0,   $
we obtain
\begin{equation}
    \label{sol_hopf}
    \alpha_c = \frac{rm_3(m_2 + n_2)}{\gamma_1n_2zK},
\end{equation}
where
\begin{align*}
    z &= -\frac{h_1}{2} + \sqrt{\frac{h_1^2}{4} + h_2}\,, \\
    h_1 &= m_2 + m_3 + n_2 + \frac{m_3(m_2 + n_2)}{m_2 + m_3 + n_2}, \\
    h_2 &= \frac{m_3(m_2 + n_2)(r - m_1 - n_1)}{m_2 + m_3 + n_2}.
\end{align*}

Considering condition (b) of the \textit{Liu's criterion}, we have
\begin{equation*}
    \left(\frac{d\Delta}{d\alpha}\right)_{\alpha = \alpha_c} = -\frac{B_1}{\alpha_c^2} - \frac{2B_2}{\alpha_c^3} < 0,
\end{equation*}
where
\begin{align*}
    B_1 &= \frac{rm_3(m_2 + n_2 + m_3)^2(m_2 + n_2)+rm_3(m_2+n_2)^2}{\gamma_1n_2K}, \\
    B_2 &= \frac{r^2(m_2 + n_2 + m_3)(m_2 + n_2)^2+m_3^2}{\gamma_1^2n_2^2K^2}.
\end{align*}

Hence, according to Liu's criterion, a simple Hopf bifurcation occurs at $\alpha = \alpha_c$, that is, $A_5 = 1$.

\section{Seasonal dynamics of the sugarcane borer-parasitoid agroecosystem}

Several environmental parameters (such as air temperature, air humidity, rainfall dispersion, among others) fluctuate periodically affecting an ecological system dynamics.
Thus, they can be represented as periodic-time functions. In this section, the intrinsic growth rate $r$ in system \eqref{model} is considered as a sinusoidal function representing these seasonal perturbations. The parameter $r$ being defined by the following function \cite{rinaldi1993multiple, gakkhar2003chaos, gakkhar2003seasonally, altizer2006seasonality}:

\begin{equation}
    \label{sas_function}
    r(t) = r_0\left(1 + r_1\sin\left(\frac{2\pi t}{365}\right)\right),
\end{equation}
where $t$ is measured in days, so $r_0$ is the average value of $r$ over an integer number of years. The parameter $r_1$ represents the degree of seasonality, hence $r_0r_1$ is the magnitude of the perturbation in $r$.

Next, we are interested in the seasonal dynamics of the equilibrium point $E_5$, where there is coexistence of the parasitoid and pest populations. For this, we will keep the values of the parameters fixed as follows \cite{parraetal2002,rafikov2012mathematical,rafikov2014dynamical}.

\begin{align}
    \label{param_values}
    m_1 &= 0,\ m_2 = 0.03566,\ m_3 = \frac{1}{4},\ m_4 = 0.00257,\ m_5 = m_4,\ m_6 = \frac{1}{5},\ \nonumber \\ n_1 &= \frac{1}{8} ,\
    n_2= \frac{1}{9},\ n_3 = \frac{1}{50},\ n_4 = \frac{1}{16},\ r_0 = 0.19,\ \\ \beta &= 0.000009,\ \gamma_1 = 2.29,\ \gamma_2 = 40,\ K = 25000 \nonumber. 
\end{align}

Setting the parameter values as specified in \eqref{param_values}, it can be shown from \eqref{equilibrio5} that without seasonality, that is, $r_1 = 0$ in \eqref{sas_function}, depending on the value of $\alpha > 0.169 \times 10^{-4}$, the attractor in the phase space can be an equilibrium point or a limit cycle, namely:
\begin{itemize}
    \item For $0.169 \times 10^{-4} < \alpha < \alpha_c$, where $\alpha_c = 0.9135 \times 10^{-4}$, the corresponding attractor is the coexistence equilibrium point $E_5$, whose components depend on the value of $\alpha$, as specified in \eqref{eq5coor};
    \item For $\alpha > \alpha_c$, that is, beyond the Hopf bifurcation, the attractor is a period one limit-cycle and, the amplitude of this limit cycle in the phase space increases with increasing the value of $\alpha$.
\end{itemize}

The addition of seasonality to the model,  through the population growth rate according to \eqref{sas_function}, induces a destabilizing effect and may even trigger chaotic behavior. This destabilizing effect will be confirmed further by computer simulations with the parameter values fixed according to \eqref{param_values}.

\begin{figure}[!h]
\centering
  \includegraphics[width=17cm, height=6cm]{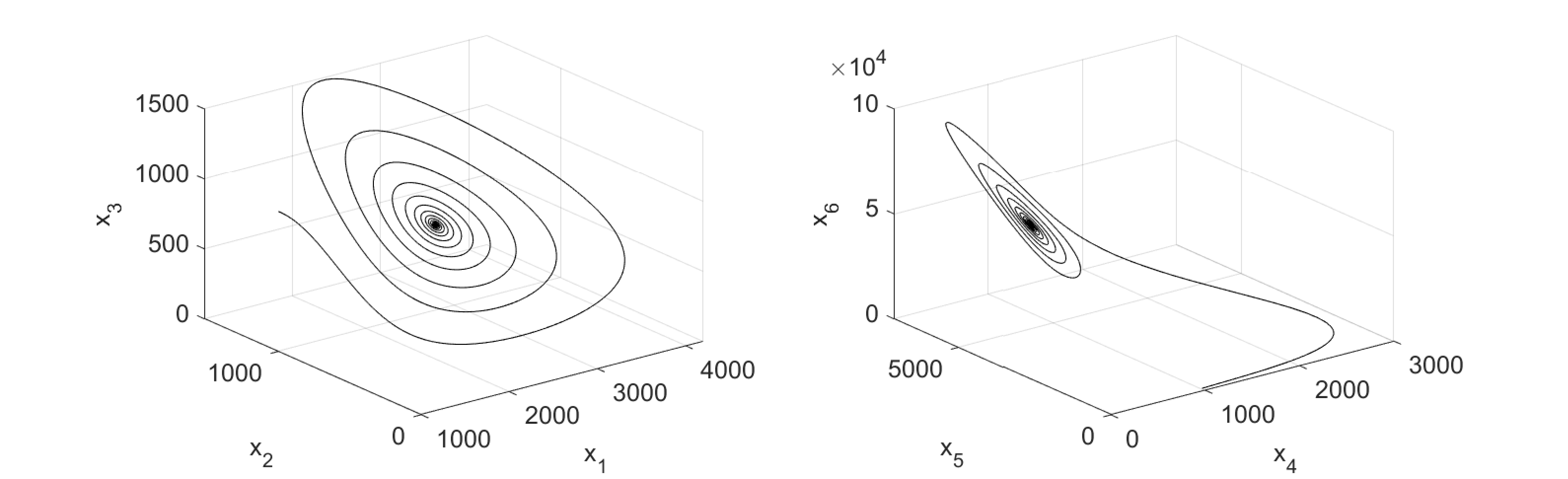}
\caption{The equilibrium $E_5$, reached by the populations of system \eqref{model}, without seasonality, for the parameters fixed in \eqref{param_values} and $\alpha = 0.6 \times 10^{-4} < \alpha_c$.}
\label{fig44}       
\end{figure}

\begin{figure}[!h]
    \centering
        \includegraphics[height=5cm,width=10cm]{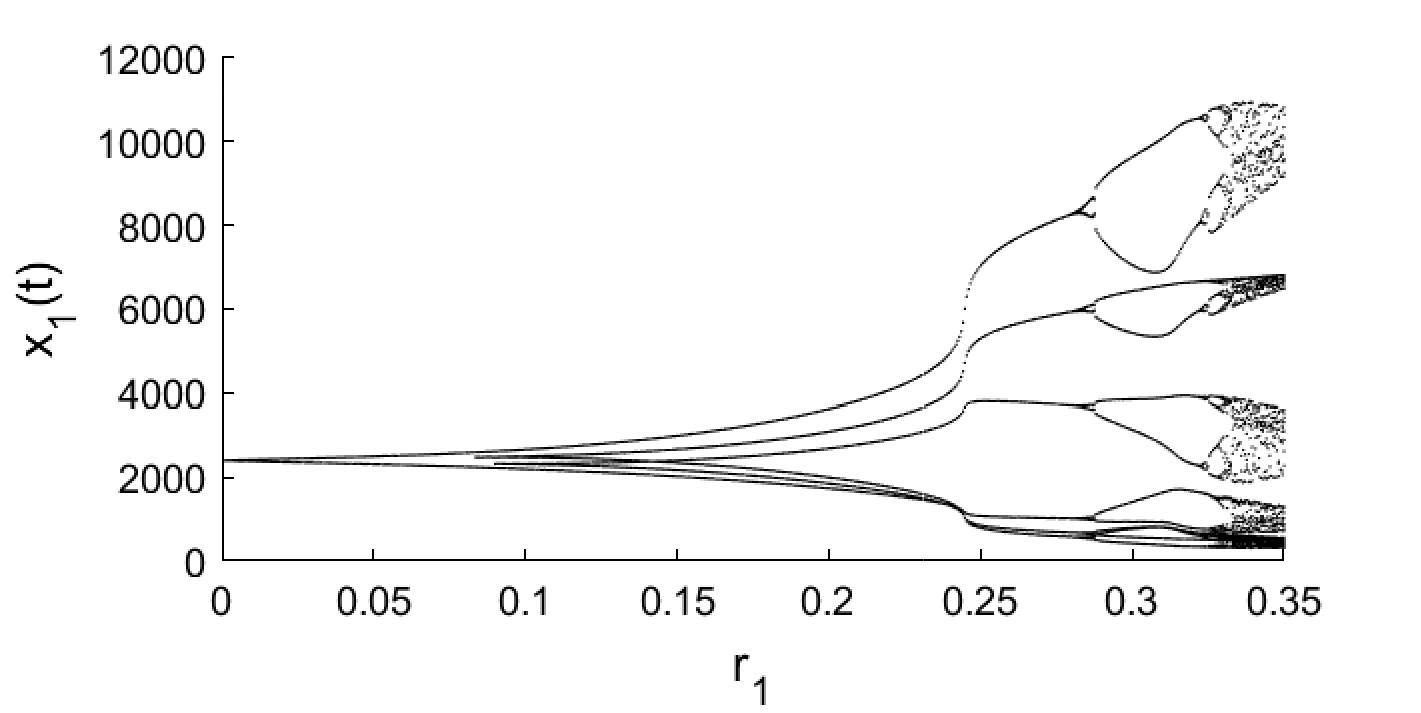}
    \caption{Keeping the parameter values fixed in \eqref{param_values} and $\alpha = 0.6 \times 10^{-4}$, the bifurcation diagram of $x_1(t)$ for $0 \leq r_1 \leq 0.35$.} 
    \label{e5_sas_periodic}
\end{figure}

\begin{figure}[!h]
\centering
    \includegraphics[width=17cm, height=6cm]{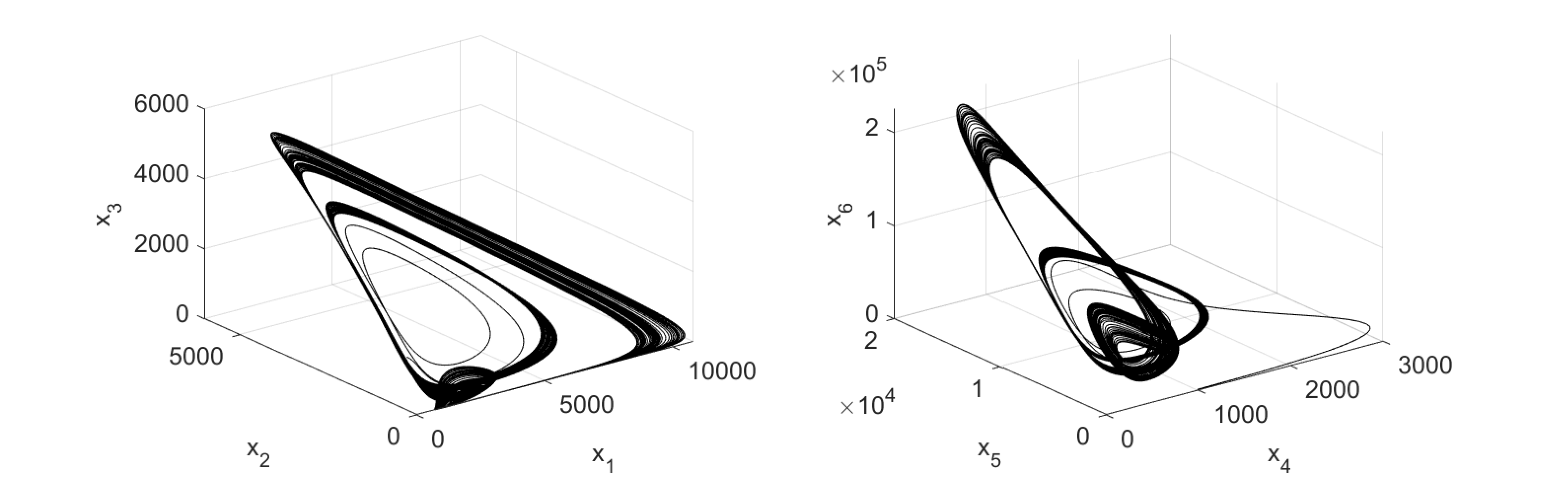}
\caption{Keeping the parameter values fixed according to \eqref{param_values} and $\alpha = 0.6 \times 10^{-4}$, the  projections of the
6D strange attractor in the phase space, corresponding to $r_1 = 0.35$, in the subspaces $x_1x_2x_3$ and $x_4x_5x_6$.}
\label{fig55}       
\end{figure}

First, we investigate the effect of seasonality for a value of $\alpha$ when the attractor is the equilibrium point $E_5$ in the 6D phase space, as shown in Fig. \ref{fig44} for $\alpha = 0.6 \times 10^{-4}$. Considering this value of $\alpha$, the bifurcation diagram for $0 \leq r_1 \leq 0.35$ is shown in Fig. \ref{e5_sas_periodic}, where it can be immediately verified that the value of $x_1(t)$ for $r_1 = 0$ is, as expected, the same value of this component of the equilibrium point obtained without seasonality, given by \eqref{eq5coor}, and shown in Fig. \ref{fig44}. Increasing the value of $r_1$, there is a periodic solution followed by a period doubling sequence, which constitutes a route to chaos that occurs for $r_1 > 0.33$. The projections of the 6D strange attractor, in the phase space, with $r_1 = 0.35$, are plotted in Fig. \ref{fig55}. Thus, the seasonality destabilizes the $E_5$ equilibrium, changing the attractor from equilibrium point to limit cycle, and then to more complex dynamics as $r_1$ increases.

\begin{figure}[!h]
\centering
    \includegraphics[width=13cm, height=6cm]{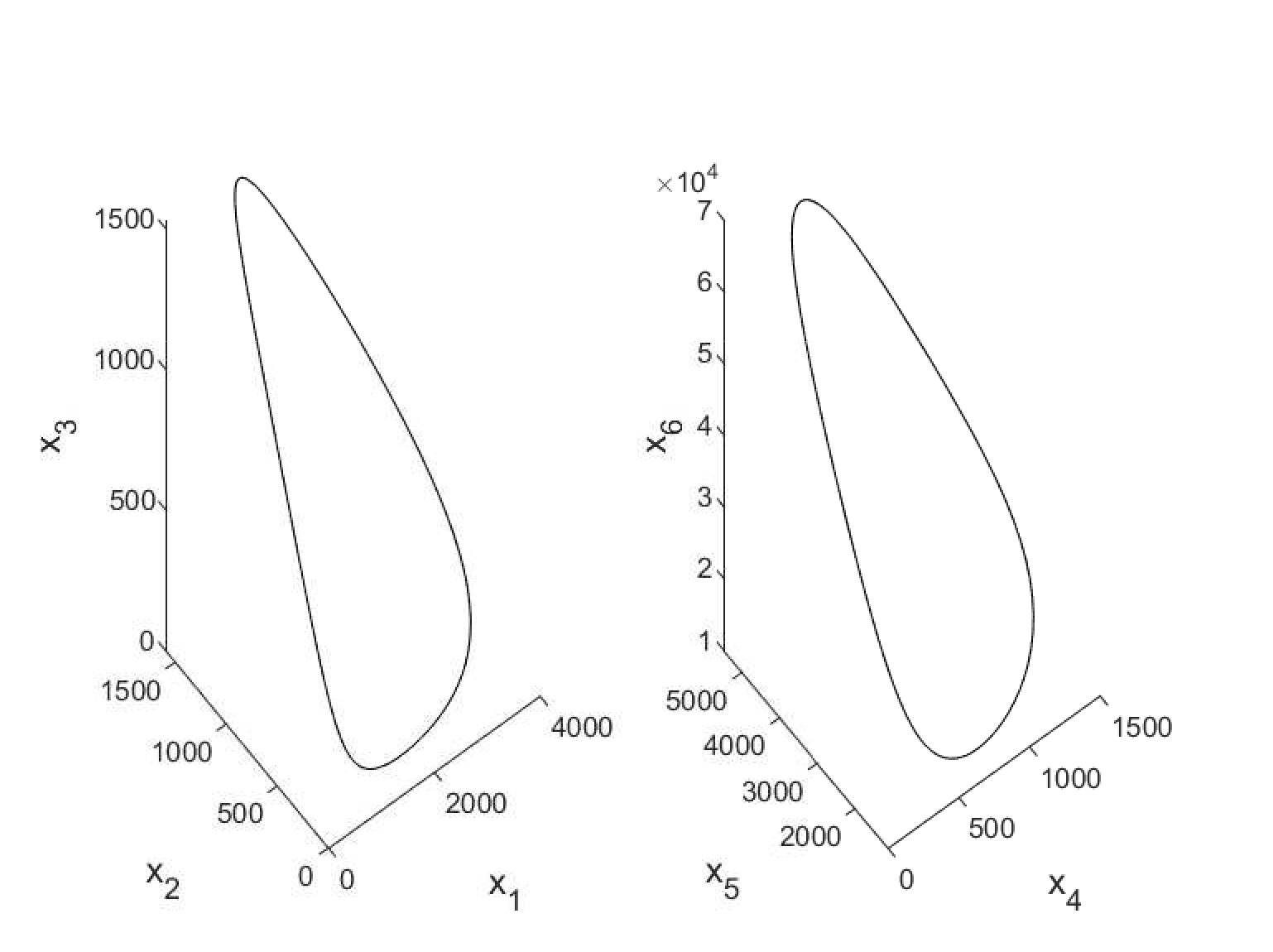}
\caption{The projection of the 6D limit cycle reached by the populations of system \eqref{model}, without seasonality, for the parameters fixed in \eqref{param_values} and $\alpha = 1 \times 10^{-4} > \alpha_c$, in the subspaces $x_1x_2x_3$ and $x_4x_5x_6$.}
\label{newfig4}       
\end{figure}

\begin{figure}[!h]
    \centering
        \includegraphics[height=5cm,width=10cm]{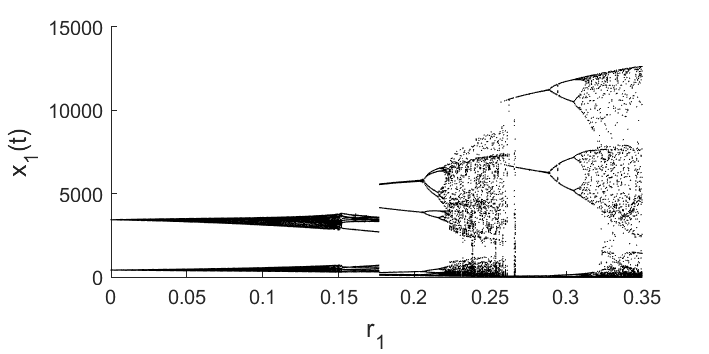}
    \caption{Keeping the parameter values fixed in \eqref{param_values} and $\alpha = 1 \times 10^{-4}$, the bifurcation diagram of $x_1(t)$ for $0 \leq r_1 \leq 0.35$.} 
    \label{limcyc_sas}
\end{figure}

\begin{figure}[!h]
\centering
  \includegraphics[width=12cm, height=5cm]{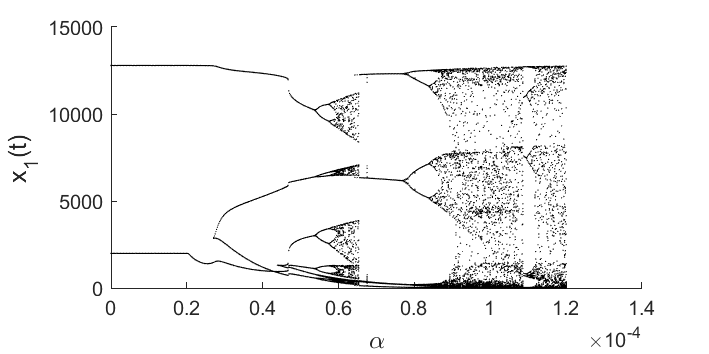}
\caption{Keeping the parameter values fixed in \eqref{param_values} and $\alpha = 1 \times 10^{-4}$, the  projections of the
6D strange attractor in the phase space, corresponding to $r_1 = 0.25$, in the subspaces $x_1x_2x_3$ and $x_4x_5x_6$.}
\label{fig77}       
\end{figure}

\begin{figure}[!h]
\centering
  \includegraphics[width=15cm, height=5cm]{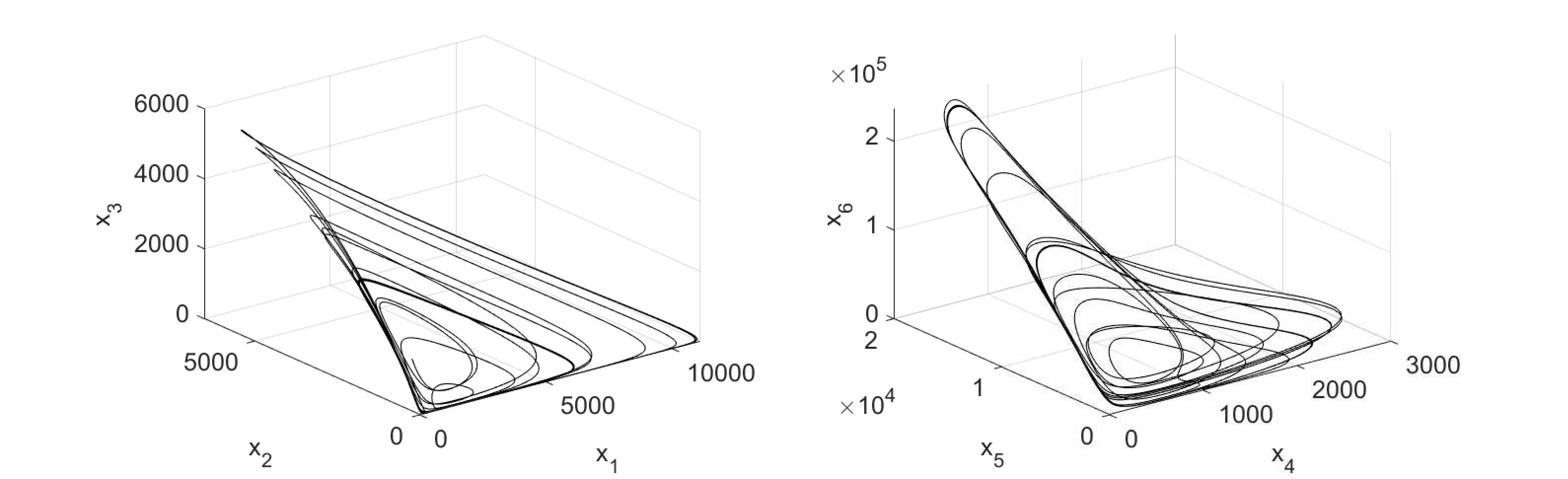}
\caption{Keeping the parameter values fixed in \eqref{param_values} and $\alpha = 1 \times 10^{-4}$, the  projections of the
6D periodic attractor in the phase space, corresponding to $r_1 = 0.28$, in the subspaces $x_1x_2x_3$ and $x_4x_5x_6$.}
\label{fig89}       
\end{figure}

Now, the range $\alpha > \alpha_c$ is considered, with the system having a period one limit cycle in the phase space, whose amplitude increases as $\alpha$ increases. We consider $\alpha = 1 \times 10^{-4}$ in Fig. \ref{newfig4}, and the maximum and minimum values of $x_1$ as the same that as we identify for $r_1 = 0$ in the bifurcation diagram of $x_1(t)$ plotted with seasonality in Fig \ref{limcyc_sas} for $0 \leq r_1 \leq 0.35$. Increasing the value of $r_1$, a period doubling sequence can be noted, which constitutes a route to chaos that occurs for $r_1 > 0.2$. The projections
of the 6D strange attractor in the phase space, corresponding to $r_1 = 0.25$, are plotted in Figure \ref{fig77}. Increasing even more the value of $r_1$, periodic attractors emerge, as the one plotted in Fig. \ref{fig89}, corresponding to $r_1 = 0.28$, and similar periodic attractors occur for $0.26 \leq r_1 \leq 0.3$. Therefore, the seasonality destabilizes the attractor from period one limit cycle, changing the behavior of our system to more complex dynamics as $r_1$ increases.

Regarding to the value of $r_1$ at which chaos is observed to occur, the comparison of bifurcation diagrams in Figs. \ref{e5_sas_periodic} and \ref{limcyc_sas} show that if $\alpha > \alpha_c$, chaos occurs at a lower value of $r_1$ than if $\alpha < \alpha_c$.

\begin{figure}[!h]
    \centering
        \includegraphics[height=8cm,width=11cm]{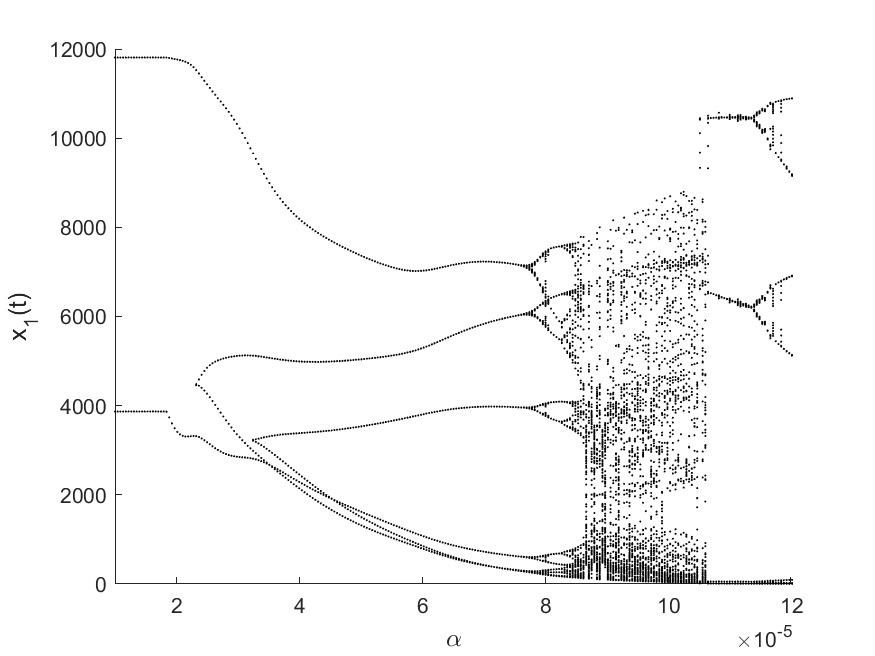}
    \caption{Bifurcation diagram of $x_1(t)$ for $0.169 \times 10^{-4} \leq \alpha \leq 1.2 \times 10^{-4}$, and keeping the parameter values fixed according to \ref{param_values}, for $r_1 = 0.25$.}
    \label{bifurcalpha}
\end{figure}

\begin{figure}[!h]
    \centering
        \includegraphics[height=8cm,width=11cm]{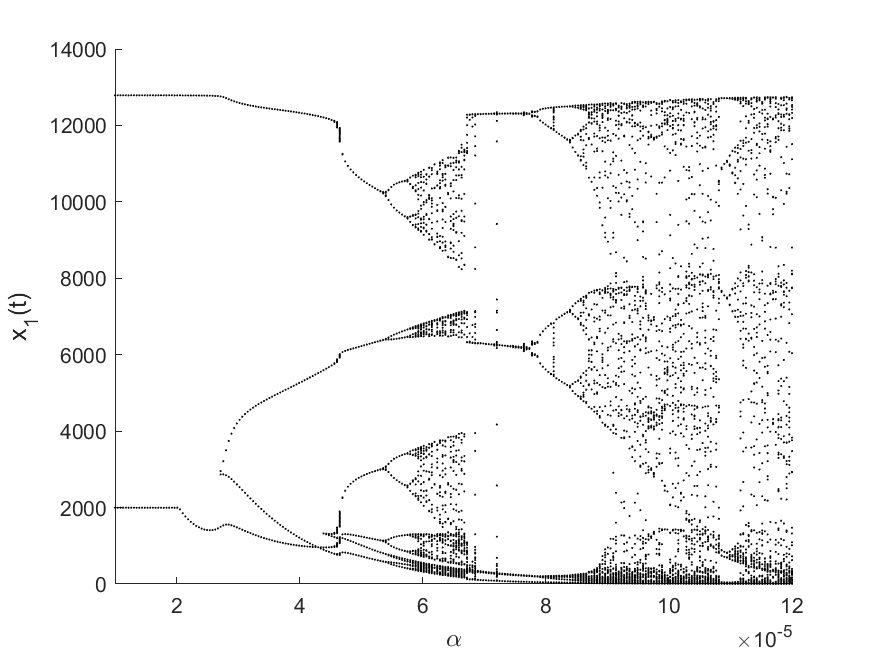}
    \caption{Bifurcation diagram of $x_1(t)$ for $0.169 \times 10^{-4} \leq \alpha \leq 1.2 \times 10^{-4}$, and keeping the parameter values fixed according to \ref{param_values}, for $r_1 = 0.35$.}
    \label{bifurcalph2}
\end{figure}

Additionally, we can investigate the effect of seasonality for fixed values of its degree $r_1$, while varying the parameter $\alpha$. For that, bifurcation diagrams of $x_1(t)$ are plotted in Figures \ref{bifurcalpha} and \ref{bifurcalph2}, keeping the parameter values fixed in \eqref{param_values} for $0.169 \times 10^{-4} \leq \alpha \leq 1.2 \times 10^{-4}$, and setting values for the parameter $r_1$. The diagram presents in Fig. \ref{bifurcalpha} corresponds to $r_1 = 0.25$, whose strange attractor considering $\alpha = 1 \times 10^{-4}$ was visualized in Fig. \ref{fig77}, while the diagram in Fig. \ref{bifurcalph2} the degree of seasonality is $r_1 = 0.35$, whose strange attractor for $\alpha = 0.6 \times 10^{-4}$ was visualized in Fig. \ref{fig55}. Comparing these two bifurcation diagrams, we conclude that the higher value of $r_1$, the lower the value of $\alpha$ at which chaos is established, as it occurs at $\alpha = 0.8 \times 10{-4}$ for $r_1 = 0.25$ and at $\alpha = 0.6 \times 10{-4}$ for $r_1 = 0.35$.


\section{Conclusions}

According to \cite{white2020seasonality}, seasonality is a significant feature in ecological systems driven by periodic climatic conditions, but it is often not explicitly included in either empirical or theoretical studies.
Therefore this article is an effort toward the integration of the complex dynamics involving such seasonal influences into the 
sugarcane borer agroecosystem and its parasitoids. 
In this work, we have proposed a novel, six-dimensional continuous-time dynamical system, modeling interactions between the sugarcane borer and its egg and larval parasitoid. On the analytical side, five equilibrium states and conditions for their local stability were found out. Moreover, the Lyapunov function stability analysis ensured the global asymptotical stability of the equilibrium state in which all considered populations coexist. Then, the occurrence of a Hopf bifurcation was investigated applying Liu’s theorem. 
Numerical simulations revealed the chaotic behavior of the system with seasonality. These results show how seasonality changes considerably the agroecosystem dynamics leading an asymptotically stable system, as shown in Fig \ref{fig44} to the period-doubling and subsequently to a chaotic attractor shown in Fig \ref{e5_sas_periodic}. For a real system, this means sudden changes can occur in populations that without seasonality could coexist in an equilibrium, in the presence of seasonal conditions. Moreover, when populations exhibit periodic oscillations, as shown in Fig. \ref{newfig4}, the introduction of seasonality can transform these oscillations into chaos, even for smaller values of $r_1$ then in the case of Fig. \ref{limcyc_sas}. Finally, bifurcation diagrams of the maximum and minimum population values, in the presence of seasonal influences,  show that an increase in the parasitism coefficient $\alpha$ can lead the stabilized system to a chaotic regime, as shown by Figures \ref{bifurcalpha} and \ref{bifurcalph2}.

The present results help understand the dynamics of the six-dimensional agro-ecological system with the seasonal forcing. Using these results the biological control strategies can be investigated in future research.


\bibliography{bibliography}

\end{document}